\documentclass{amsart}
\usepackage{amsmath}
\usepackage{amssymb}
\usepackage{graphicx}

\newtheorem{theorem}{Theorem}[section]
\newtheorem{lemma}[theorem]{Lemma}

\theoremstyle{definition}
\newtheorem{definition}[theorem]{Definition}

\theoremstyle{remark}

\def\co{\colon\thinspace}

\numberwithin{equation}{section}

\begin{document}

\title[Bordered $\widehat{HF}$ and Lefschetz fibrations with corners]{Bordered Floer Homology and Lefschetz fibrations with corners}


\author{Tova Brown}
\address{}
\curraddr{}
\email{}
\thanks{The author wishes to thank their thesis adviser, Denis Auroux; Matt Hedden, for instruction in Heegaard Floer theory and other topics; and Robert Lipshitz, for many helpful conversations.}



\date{}

\begin{abstract}

Lipshitz, Ozsv\'ath and Thurston defined a bordered Heegaard Floer invariant $\widehat{CFDA}$ for 3-manifolds with two boundary components, including mapping cylinders for surface diffeomorphisms.  We define a related invariant for certain 4-dimensional cobordisms with corners, by associating a morphism $F \co \widehat{CFDA}(\phi_1) \to \widehat{CFDA}(\phi_2)$ to each such cobordism between two mapping cylinders $\phi_1$ and $\phi_2$.  Like the Osv\'ath-Szab\'o invariants of cobordisms between closed 3-manifolds, this morphism arises from counting holomorphic triangles on Heegaard triples.  We demonstrate that the homotopy class of the morphism $F$ only depends on the symplectic structure of the cobordism in question.  

\end{abstract}

\maketitle

\section{Introduction}

Heegaard Floer theory is a set of invariants for closed, connected 3-manifolds and cobordisms between them, with a related invariant for closed 4-manifolds \cite{OS1, OS2}.  Together these invariants form a $3+1$ dimensional topological quantum field theory (TQFT), meaning a functor from the cobordism category of 3-manifolds to, in this case, the category of graded abelian groups.  \\

The construction of Heegaard Floer homology involves counting holomorphic curves associated to Heegaard diagrams of 3-manifolds.  Specifically, given a 3-manifold $Y$ with a genus $g$ Heegaard diagram $\{H, \alpha_i, \beta_i\}$, the invariant $\widehat{HF}(Y)$ is defined as the homology of a chain complex generated by g-tuples of intersection points between the $\alpha$ and $\beta$ curves.  In Lipshitz' reformulation \cite{Lipshitz}, the differential arises from counts of rigid holomorphic curves in the symplectic manifold $H \times [0,1] \times \mathbb{R}$, with boundaries mapping to the Lagrangian submanifolds $\beta_i \times \{0\} \times \mathbb{R}$ and $\alpha_i \times \{1\} \times \mathbb{R}$.  The maps associated to cobordisms arise from a similar construction, which uses Heegaard triples to represent certain elementary cobordisms \cite{OS2}.  \\

In 2008, Lipshitz, Ozsv\'ath and Thurston \cite{LOT1} developed bordered Heegaard Floer homology, which generalizes $\widehat{HF}$ to parametrized Riemann surfaces and to bordered 3-manifolds, meaning 3-manifolds with parametrized boundary.  Given two such 3-manifolds $Y_1$ and $Y_2$, if the surfaces $\partial(Y_1)$ and $\partial(Y_2)$ have compatible parametrizations, then the bordered Heegaard Floer invariants for $Y_1$ and $Y_2$ may be combined to obtain $\widehat{HF}(Y)$, where $Y$ is the 3-manifold defined by identifying the boundaries of $Y_1$ and $Y_2$.  \\

Specifically, to a parametrized surface $S$, there is an associated differential graded algebra $A_S$. If $\partial(Y_1)$ is identified with $S$ and $\partial(Y_2)$ with $-S$, then the bordered invariant for $Y_1$ is a right $A_{\infty}$ module $\widehat{CFA}(Y_1)$ over $A_S$, while the invariant for $Y_2$ is a left differential graded module with an additional ``type D'' structure over $A_S$, called $\widehat{CFD}(Y_2)$.  Lipshitz, Ozsv\'ath and Thurston define the tensor product $\widehat{CFA}(Y_1) \boxtimes \widehat{CFD}(Y_2)$, which is a simple model for the $A_{\infty}$ tensor product.  They then demonstrate that the resulting chain complex is quasi-isomorphic to the closed invariant $\widehat{HF}(Y)$.  \\

Given such a decomposition of a closed 3-manifold $Y = Y_1 \cup Y_2$, we may represent $Y$ by a Heegaard diagram $H = H_1 \cup H_2$, where $H_1$ and $H_2$ are subsurfaces of $H$ with disjoint interiors, each $\beta$ curve is contained entirely in either $H_1$ or $H_2$, and $Y_i$ is the union of all gradient flow lines of the Morse function that pass through $H_i$, for each $i$.  The marked surfaces $H_1$ and $H_2$ are called bordered Heegaard diagrams for $Y_1$ and $Y_2$, and they contain the data needed to define $\widehat{CFA}(Y_1)$ and $\widehat{CFD}(Y_2)$, respectively.  \\

In each case, the generators are the tuples of intersection points of the $\alpha$ and $\beta$ curves in $H_i$ which extend to generators of $\widehat{CF}(Y)$, while the differential and products involve counting rigid holomorphic curves.  However, in order to rebuild the closed invariant from these pieces, the algebra $A_S$ and the modules $\widehat{CFA}(Y_1)$ and $\widehat{CFD}(Y_2)$ must encode information about how such curves interact with the boundary $H_1 \cap H_2$.  To accomplish this, the generators of $A_S$ are ``strand diagrams'' representing ways that rigid holomorphic curves may intersect $H_1 \cap H_2$, while the relations in $A_S$ represent ways that the ends of one-dimensional moduli spaces of holomorphic curves may behave near this boundary.  \\

In the $A_{\infty}$ module $\widehat{CFA}$, the products record the behavior of holomorphic curves that hit the boundary in certain prescribed ways, with rigid curves that intersect the boundary more times contributing to higher products.  The type $D$ structure on $\widehat{CFD}$ consists of a differential and an identification between $\widehat{CFD}$ and $A_S \otimes \chi$, where $\chi$ is the $\mathbb{Z}/2\mathbb{Z}$ vector space whose generators are the same as those of $\widehat{CFD}$, with this data satisfying certain properties.  \\

Lipshitz, Ozsv\'ath and Thurston also defined a bordered invariant for cobordisms between parametrized surfaces \cite{LOT2}.  This is a bimodule, called $\widehat{CFDA}$, which incorporates both the type D structure and the $A_{\infty}$ structures of the modules $\widehat{CFD}$ and $\widehat{CFA}$.  Bimodules with this structure are called type $DA$ bimodules.  \\

The bimodule $\widehat{CFDA}$ is defined for 3-dimensional cobordisms in general, but in particular we may consider mapping cylinders of surface diffeomorphisms, meaning 3-manifolds diffeomorphic to a product $\Sigma \times [0,1]$ with the boundary components parametrized, and with a marked, framed section over $[0,1]$ which allows us to compare the two parametrizations.  This yields a functor from the mapping class groupoid to the category of differential graded algebras, with morphisms given by type $DA$ bimodules.  \\

We may construct a 2-category from the mapping class groupoid by taking certain Lefschetz fibrations over rectangles as 2-morphisms.  The main result of this paper is that these cobordisms induce type $DA$ maps between the $\widehat{CFDA}$ invariants of mapping cylinders, and that this data forms a 2-functor.  \\

Specifically, the 2-morphisms we use are ``cornered Lefschetz fibrations," or CLF's.  A CLF is a Lefschetz fibration over a rectangle with certain markings on its fibers.   The left and right edges are identified with $\Sigma_i \times I$ for some parametrized surfaces $\Sigma_1$ and $\Sigma_2$, respectively, while the top and bottom edges are identified with mapping cylinders, so the resulting parametrizations of the corners coincide.  This Lefschetz fibration is also equipped with a marked framed section, which corresponds to the marked sections on the edges.  With this definition understood, we have the following theorem:  

\begin{theorem}
Given a cornered Lefschetz fibration between two mapping cylinders $f$ and $f'$, there is an induced type DA bimodule map from $\widehat{CFDA}(f)$ to $\widehat{CFDA}(f')$.  This map is well-defined up to chain homotopy.  
\end{theorem}

To define a cobordism map associated to a CLF with a single critical point, we first construct a bordered Heegaard triple which represents this Lefschetz fibration.  To accomplish this, we consider the vanishing cycle as a knot in the mapping cylinder identified with the bottom edge, and build a genus $2g$ bordered Heegaard diagram for this mapping cylinder subordinate to that knot, where $g$ is the genus of the fiber.  We then define an additional set of curves, obtaining a Heegaard triple which represents the cobordism induced by the appropriate surgery.  \\

The cobordism map is defined by counting rigid holomorphic triangles associated with this Heegaard triple.  The higher maps and type $D$ structure maps encode the ways that these triangles interact with the right and left boundaries of the Heegaard surface, respectively.  \\

More generally, we may associate a cobordism map to any CLF, by decomposing this Lefschetz fibration into pieces by a sequence of horizontal and vertical cuts.  Given two CLF's $W_1$ and $W_2$, with the right edge of $W_1$ and the left edge of $W_2$ equipped with compatible parametrizations, we may define their horizontal composition $W_1 \circ_h W_2$ by identifying these edges.  If $W_i$ is a cobordism between the mapping cylinders $f_i$ and $f_i'$, and we have maps $F_i: \widehat{CFDA}(f_i) \rightarrow \widehat{CFDA}(f_i')$  associated to each $W_i$, then there is an induced type $DA$ bimodule map: $$F_1 \boxtimes F_2: \widehat{CFDA}(f_1) \boxtimes \widehat{CFDA}(f_2) \rightarrow \widehat{CFDA}(f_1') \boxtimes \widehat{CFDA}(f_2').$$  

Similarly, if $V_1$ and $V_2$ are CLF's where $V_1$ is a cobordism from $f$ to $g$ and $V_2$ is a cobordism from $g$ to $h$, then we may define the vertical composition $V_1 \circ_v V_2$ by identifying the top edge of $V_1$ with the bottom edge of $V_2$.  Given maps between the appropriate bimodules associated to $V_1$ and $V_2$, we may associate the composition of these maps to the vertical composition of $V_1$ and~$V_2$.  \\

To prove that the homotopy class of maps associated to a CLF with multiple critical points does not depend on the decomposition, we will show that horizontal decompositions may be altered to form vertical decompositions, and vice versa.  This flexibility allows us to show that Hurwitz moves do not change the class of the map, and also allows us to rearrange a description of a given Lefschetz fibration in order to facilitate the calculation of the invariant.  

\section{Type $DA$ bimodules}

In this section we will review the concepts of $A_{\infty}$ modules, type $D$ modules, and type $DA$ bimodules, working exclusively over $\mathbb{Z}/2\mathbb{Z}$.  This material is covered in greater detail in Section 2 of \cite{LOT1}, and Section 2 of \cite{LOT2}.  

\subsection{$A_{\infty}$ modules}

Given a differential graded algebra $A$ with differential $\delta$, a right $A_{\infty}$ module over $A$ is vector space $M$ over $\mathbb{Z}/2\mathbb{Z}$, with a differential $m_1: M \rightarrow M$ and products $m_{n+1}: M \otimes A^{\otimes n} \rightarrow M$ for $n \geq 1$, satisfying the property:

\begin{align*}
0 &= \sum_{i+j=n} m_{j+1}(m_{i+1}(x, a_1, \ldots a_i), a_{i+1}, \ldots a_n) \\
&+ \sum_{1 \leq k \leq n} m_{n+1}(x, a_1, \ldots a_{k-1}, \delta(a_k), a_{k+1}, \ldots a_n) \\
&+ \sum_{1 \leq k < n} m_n(x, a_1, \ldots a_{k-1}, a_k a_{k+1}, a_{k+2}, \ldots a_n),
\end{align*}

for each $n \geq 0$.  Note that by taking $n=0$ or $n=1$ we obtain the familiar rules $m_1 \circ m_1 = 0$ and $m_1(m_2(x, a_1)) = m_2(m_1(x), a_1) + m_2(x, \delta(a_1))$.  By taking $n = 2$ we see that, while the product $m_2$ need not be associative, it does associate up to a chain homotopy given by $m_3$.  In general, these properties ensure that each $m_n$ resolves the failures of associativity which arise in the products $m_i$, for $i < n$.  \\

We may also define a right $A_{\infty}$ module over an $A_{\infty}$ algebra $A$.  Here, $A$ is a vector space over $\mathbb{Z}/2\mathbb{Z}$, equipped with a differential $\mu_1$ and products $\mu_n: A^{\otimes n} \rightarrow A$ for $n \geq 2$.  For each $n$, these operations satisfy:

$$0 = \sum_{i+j+1=n, 1 \leq k \leq j+1} \mu_{j+1}(a_1, \ldots a_{k-1}, \mu_{i+1}(a_k, \ldots a_{k+i}), a_{k+i+1}, \ldots a_n).$$  

Given such an $A_{\infty}$ algebra, a right $A_{\infty}$ module $M$ over $A$ is a $\mathbb{Z}/2\mathbb{Z}$ vector space, equipped with a differential $m_1$ and products $m_{n+1}:M \otimes A^{\otimes n} \rightarrow M$, satisfying the properties:  

\begin{align*}
0 &= \sum_{i+j=n} m_{j+1}(m_{i+1}(x, a_1, \ldots a_i), a_{i+1}, \ldots a_n) \\
&+ \sum_{i+j=n, 1 \leq k \leq j} m_{j+1}(x, a_1, \ldots a_{k-1}, \mu_{i+1}(a_k, \ldots a_{i+k}), a_{i+k+1}, \ldots a_n).
\end{align*}

While this paper only makes use of honest differential graded algebras and related structures, the definitions of these objects often become clearer when viewed in a more general context.  Later in this section we will define type $DA$ bimodules over $A_{\infty}$ algebras, noting that these definitions reduce to those of analogous structures over DGA's.  

\subsection{Type $D$ modules}

Lipshitz, Ozsv\'ath and Thurston define a structure called a type D module; see Definition 2.12 in \cite{LOT1} and Definition 2.2.20 in \cite{LOT2}.  This is a left module $N$ over a differential graded algebra $A$, with a differential $\delta$ and an identification $D \co N \rightarrow A \otimes \chi$, where $\chi$ is the $\mathbb{Z}/2\mathbb{Z}$ vector space generated by the generators for $N$.  \\

The existence of this identification allows us to study the behavior of the differential in greater detail.  The differential satisfies $\delta^2 = 0$, and so repeating this map does not yield information.  However, we may consider instead the map $D_1 \co N \rightarrow A \otimes \chi,$ defined by $D_1 = D \circ \delta$.  \\

Given an element $x \in N$, the element $D_1(x)$ is of the form $\sum a_i \otimes y_i$, for some $a_i \in A$ and $y_i \in \chi$.  While $\delta(x)$ is a cycle the elements $y_i$ may not be, and so we may consider the element of $A \otimes A \otimes \chi$ given by $\sum a_i \otimes D \circ \delta(y_i)$.  We may repeat the above process an arbitrary number of times, obtaining an element $D_n(x) \in A^{\otimes n} \otimes \chi$ for each $x \in N$, defined recursively by $D_n = (\mathbb{I}_{A^{\otimes n-1}} \otimes D_1) \circ D_{n-1}$.  \\

The maps $D_n$ satisfy a set of properties ensuring that $\delta^2 = 0$.  In an $A_{\infty}$ module the products $m_n$ satisfy properties ensuring associativity up to homotopy, and so the $D_n$ and $m_n$ play similar roles.  In both cases we wish to examine the behavior of a mechanism which resolves an ambiguity, and the $m_n$ and $D_n$ describe the specifics of that mechanism.  \\

Given a differential graded algebra $A$, a right $A_{\infty}$ module $M$ over $A$, and a type D module $N \simeq A \otimes \chi$, Lipshitz, Ozsvath and Thurston define a differential $d$ on $M \otimes \chi$ (see section 2.3 of \cite{LOT1}), obtaining a chain complex $M \boxtimes N$.  This differential arises from the relationships between the products $m_n$ and the maps $D_n$ for particular generators, namely:  

$$d = \sum (m_{n+1} \otimes \mathbb{I}_{\chi}) \circ (\mathbb{I}_M \otimes D_n).  $$

\subsection{Type $DA$ bimodules}

Given $A_{\infty}$ algebras $A_1$ and $A_2$, a type $DA$ bimodule over $A_1$ and $A_2$ is an object which behaves as an $A_{\infty}$ module over $A_2$ and a type $D$ module over $A_1$ (see section 2.2.4 of \cite{LOT2}).  It is a $\mathbb{Z}/2\mathbb{Z}$ vector space $M$ with a differential $m_1$ and right products $m_{n+1} : M \otimes A_2^{\otimes n} \rightarrow M$, satisfying the usual $A_{\infty}$ properties.  $M$ is equipped with an identification $D : M \rightarrow A_1 \otimes \chi$, where $\chi$ is a $\mathbb{Z}/2\mathbb{Z}$ vector space, and this endows $M$ with a left product by elements of $A_1$.  \\

Just as the function $D$ associated with a type $D$ module allows us to examine the behavior of the differential, the analogous map associated with a type $DA$ bimodule $M$ allows us to study the products $m_n$ in greater detail.  Repeated products in a type $DA$ bimodule are constrained by the $A_{\infty}$ relations, since, for example, the term $m_{j+1}(m_{i+1}(x, a_1, \ldots a_i), a_{i+1}, \ldots a_{i+j})$ appears in the product $m_1(m_{i+j+1}(x, a_1, \ldots a_{i+j}))$.  However, if $(b, y)$ is a term in $D(m_{i+1}(x, a_1, \ldots a_i))$, then the product $m_{j+1}(y, a_{i+1}, \ldots a_{i+j})$ is not so constrained.  \\

More generally, for each $n \geq 1$ we have a map:  

$$D_n: \chi \otimes T^*A_2 \rightarrow A_1^{\otimes n} \otimes \chi.$$

These maps are defined recursively, by:  

$$D_1(x, a_1, \dots a_i) = D(m_{i+1}(x, a_1, \ldots a_i)),$$

$$D_n(x, a_1, \dots a_i) = \sum_{0 \leq j \leq i} (\mathbb{I}_{T*A_1} \otimes D_1)(D_{n-1}(x, a_1, \ldots a_j), a_{j+1}, \ldots a_i).$$

Since the products on $M$ satisfy the $A_{\infty}$ relations, the maps $D_n$ satisfy the following property:  

$$\sum_{n \geq 1}(\mu_n \otimes \mathbb{I}_{\chi}) \circ D_n + D_1 \circ (\mathbb{I}_{\chi} \otimes m) = 0,$$

where $m$ is the differential on $T^*A_2$.  In this sense, the $D_n$ maps provide information about how the $A_{\infty}$ relations are satisfied, information which is essential when taking tensor products of type $DA$ bimodules.  

\subsection{Morphisms and chain homotopies}

Given two right $A_{\infty}$ modules $M$ and $M'$ over an $A_{\infty}$ algebra $A$, a morphism from $M$ to $M'$ consists, in part, of a chain map $F_1: M \rightarrow M'$.  If $F_1$ were a morphism of modules we would require that it preserve the product $m_2$, but in this case it need only preserve this product up to homotopy.  This means that, for each $x \in M$ and $a \in A$, there is an element $F_2(x, a) \in M'$ with:

$$m_1(F_2(x,a)) + F_2(m_1(x), a) + F_2(x, \mu_1(a)) = F_1(m_2(x,a)) + m_2(F_1(x),a).$$  

Since $M$ and $M'$ are equipped with higher products as well, a morphism between them must also preserve these products up to homotopy, and so must be equipped with a specified element $F_{n+1}(x, a_1, \ldots a_n)$ for each $(x, a_1, \ldots a_n) \in M \otimes A^{\otimes n}$ to resolve these ambiguities.  Furthermore, these higher maps $F_n$ introduce their own ambiguities which must also be resolved.  Thus we require the $F_n$ to satisfy the following property:  

\begin{align*}
0 &= \sum_{i+j = n} F_{i+1}(m_{j+1}(x, a_1, \ldots a_j), a_{j+1}, \ldots a_n) \\
&+ \sum_{i+j = n, k \leq i} F_{i+1}(x, a_1, \ldots a_{k-1}, \mu_{j+1}(a_k, \ldots a_{k+j}), a_{k+j+1}, \ldots a_n) \\
&+ \sum_{i+j = n} m_{j+1}(F_{i+1}(x, a_1, \ldots a_i), a_{i+1}, \ldots a_n).
\end{align*}

A morphism of type D modules is simply a chain map of left modules, however the presence of the type $D$ maps imposes additional structure.  Consider two type D modules $N$ and $N'$ over an algebra $A$, with structure maps $D:N \rightarrow A \otimes \chi$ and $D':N' \rightarrow A \otimes \chi'$.  A map between them is given by a function $F: \chi \rightarrow N'$, which commutes with the differentials.  However, we may also consider maps of the form:  

$$F_{i, j} = (I_{A^{\otimes j+1}} \otimes D_i') \circ (I_{A^{\otimes j}} \otimes F) \circ D_j : \chi \rightarrow A^{\otimes i+j+1} \otimes \chi'.$$

Since these maps arise from the interactions between a chain map and differentials, they must satisfy the property:  

$$\sum_{i, j} (\mu_{i+j+1} \otimes I_{\chi'}) \circ F_{i, j} = 0.$$

Now let $A_1$ and $A_2$ be $A_{\infty}$ algebras, and let $M$ and $M'$ be type $DA$ bimodules, both over $A_1$ and $A_2$.  The bimodule $M$ is equipped with products $m_i$ and a type $D$ map $D: M \rightarrow A_1 \otimes \chi$, and the bimodule $M'$ has products $m_i'$ and an identification $D': M' \rightarrow A_1 \otimes \chi'$.  \\

As defined by \cite{LOT2} in Definition 2.2.39, a type $DA$ morphism from $M$ to $M'$ is a collection of maps:

$$F_{i+1}: \chi \otimes A_2^{\otimes i} \rightarrow M'.$$

These maps must satisfy properties analogous to those of an $A_{\infty}$ morphism, while the type D structures for both $M$ and $M'$ interact with the maps $F_i$ in constrained ways.  Specifically, for each $i + j +1 = n$, we may define a map $F_{i, j} : \chi \otimes T^*A_2 \rightarrow A_1^{\otimes n} \otimes \chi'$ as follows:  

\begin{align*}
&F_{i, j}(x, a_1, \ldots a_k) = \\
&\sum_{k'=0}^k  (I_{A_1^{\otimes j+1}} \otimes D_i)(I_{A_1^{\otimes j}} \otimes F_{k'+1})(D_j \otimes I_{A_2^{\otimes k'}})(x, a_1, \ldots a_k).  
\end{align*}

Then these maps must satisfy the property:

$$\sum_n \sum_{i+j+1=n} \mu_n \circ F_{i, j} + F_{0, 0} \circ m = 0,$$

Where here $m$ is the differential on $T^*A_2.$  Observe that this requirement generalizes the properties of both $A_{\infty}$ and type $D$ morphisms.  \\

There is a notion of chain homotopies between type $DA$ morphisms, and the details of this are given in \cite{LOT2} in Definition 2.2.39.  

\subsection{Type $DA$ compositions and tensor products}

Suppose that $A_1$ and $A_2$ are $A_{\infty}$ algebras, and that  $M, M',$ and $M''$ are all type $DA$ bimodules over these algebras.  Given morphisms $F: M \rightarrow M'$ and $F': M' \rightarrow M''$, Lipshitz, Ozsv\'ath and Thurston define their composition $F' \circ F$; see Definition 2.2.39 and Figure 2 of \cite{LOT2}.  This is a type $DA$ morphism, and longer compositions $F^{(n)} \circ \ldots \circ F$ are well-defined up to homotopy.  \\

Now let $A_1, A_2$ and $A_3$ be DGA's.  Let $N$ be a type $DA$ bimodule over $A_1$ and $A_2$, and let $M$ be a type $DA$ bimodule over $A_2$ and $A_3$.  Then there is a type $DA$ bimodule $N \boxtimes M$ (Definition 2.3.8 of \cite{LOT2}), which generalizes the tensor product of $A_{\infty}$ and type $D$ modules.  Namely, if $\chi_1$ and $\chi_2$ are the generating sets for $N$ and $M$, respectively, then $N \boxtimes M$ is identified with $A_1 \otimes \chi_1 \otimes \chi_2$, and equipped with products $m'_n$ given by:

$$m'_{n+1} (x \otimes y, a_1, \ldots a_n) = \sum_{i \geq 0} (m_{i+1} \otimes I_{\chi_2}) (x, D_i(y, a_1, \ldots a_n)).$$  

Here the $m_i$ are the products associated to $N$, and the $D_i$ are the structure maps associated to the type $DA$ bimodule $M$.  \\

Given type $DA$ morphisms $F:N \rightarrow N'$ and $G:M \rightarrow M'$, there is an induced type $DA$ morphism $F \boxtimes G : N \boxtimes M \rightarrow N' \boxtimes M'$.  This is defined as:

$$F \boxtimes G= (I \boxtimes G) \circ (F \boxtimes I),$$

with the morphisms $I \boxtimes G$ and $F \boxtimes I$ as defined in Figure 5 of \cite{LOT2}.  This product operation on type $DA$ morphisms is associative up to homotopy.  \\

Now suppose we have type $DA$ morphisms $F: N \rightarrow N', F': N' \rightarrow N'', G: M \rightarrow M',$ and $G': M' \rightarrow M''$, where $N, N', N''$ are type $DA$ bimodules over $A_1$ and $A_2$, and $M, M', M''$ are type $DA$ bimodules over $A_2$ and $A_3$.  Suppose furthermore that we have type $DA$ morphisms $\hat{F}: N \rightarrow N', \hat{F'} : N' \rightarrow N''$, and $\hat{G}: M \rightarrow M'$, which are homotopic to $F, F'$, and $G$, respectively.  Then we have the following results from \cite{LOT2}:  

\begin{lemma}
The morphisms $F' \circ F$ and $\hat{F'} \circ \hat{F}$ are homotopic, and the morphisms $F \boxtimes G$ and $\hat{F} \boxtimes \hat{G}$ are homotopic.  
\end{lemma}

And:

\begin{lemma}
The induced morphisms $(F' \circ F) \boxtimes (G' \circ G)$ and $(F' \boxtimes G') \circ (F \boxtimes G)$ are equivalent up to homotopy.  
\end{lemma}

\section{Bimodules and the mapping class groupoid}

\subsection{Pointed matched circles and bordered Heegaard diagrams}

\begin{figure}
\centering
\includegraphics[width=70mm]{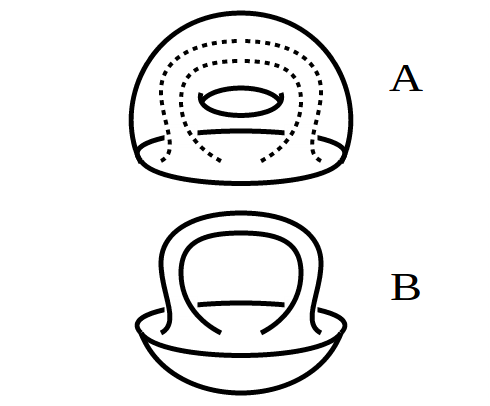}
\caption{A Heegaard decomposition of a solid torus into two handlebodies, $A$ and $B$.}  
\label{Solid_Torus}
\end{figure}

Given a genus $g$ Riemann surface $S$, a parametrization of $S$ consists of an embedded closed disk $D \subset S$, with a marked point $z$ in $\partial (D)$, along with a collection of $2g$ disjoint, properly embedded arcs $\{a_1, \ldots , a_{2g}\}$ in $\overline{S \backslash D}$, such that these arcs represent a basis for $H_1(S, D)$.  Given two parametrized Riemann surfaces $S_1$ and $S_2$, we say that their parametrizations are compatible if there is a diffeomorphism $\phi: S_1 \rightarrow S_2$ which restricts to a diffeomorphism between the marked disks, the marked points, and the collections of arcs.  If two 3-manifolds have boundary components which are parametrized in compatible ways, then such a diffeomorphism allows us to identify their boundaries in a canonical way.  \\

We may also construct a parametrized surface in the abstract, by giving a handle decomposition.  Let $c$ be an oriented circle with a marked point $z$, and with $2g$ marked pairs of points, with these points distinct from each other and from $z$.  By taking $c$ to be the boundary of a disk, we may then interpret the marked pairs as the feet of orientable 1-handles.  If no sequence of handleslides within $c \backslash \{z\}$ can bring two paired points adjacent to each other, then $c$ is called a pointed matched circle of genus $g$, and it describes a handle decomposition of a genus $g$ Riemann surface.  This surface has a canonical parametrization, in which the marked arcs are the cores of the 1-handles.  This parametrized surface is called $S_c$.  \\

A bordered 3-manifold is a 3-manifold with boundary, whose boundary components are parametrized.  Just as the Heegaard Floer invariants for closed 3-manifolds arise from Heegaard diagrams \cite{OS1}, the bordered Heegaard Floer invariants for bordered 3-manifolds arise from Heegaard diagrams with boundary \cite{LOT1}.  Let $\{H, \alpha_i, \beta_i, z\}$ be a Heegaard diagram for a manifold $Y$, and let $c$ be a separating curve on $H$ which includes the marked point $z$.  Suppose that $c$ is disjoint from the $\beta_i$, that it intersects each curve $\alpha_i$ transversely and at most twice, and that the marked pairs $\alpha_i \cap c$ make $c$ a pointed matched circle.  Then we may decompose $H$ along $c$, yielding two bordered Heegaard diagrams, $H_1$ and $H_2$.  \\

This decomposition $H = H_1 \cup H_2$ induces a decomposition of the 3-manifold $Y = Y_1 \cup Y_2$, where $Y_1$ and $Y_2$ are bordered 3-manifolds with disjoint interiors.  To see this, suppose $f: Y \rightarrow \mathbb{R}$ is a Morse function which is compatible with the Heegaard diagram $\{H, \alpha_i, \beta_i \}$.  We may then define $Y_i$ to be the closure of the union of all flow lines that pass through $H_i$, for each $i$.  \\

Let $S \subset Y$ be the closure of the union of all flow lines that pass through the curve $c$.  Then $S$ is a Riemann surface, with orientation induced by its inclusion as $\partial(Y_1)$, and this construction equips $S$ with a handle decomposition.  To see this, suppose that $f$ is self-indexing with the Heegaard surface given by $f^{-1}(\{1.5\})$.  Then $f^{-1}([0, 1.5]) \cap S$ is a closed disk with boundary $c$, and thus with the marked point $z$ on its boundary.  For each pair $\alpha_i \cap c$ the closure of its stable manifold is an arc, which we may identify with the core of a 1-handle.  The resulting parametrization is compatible with that induced by the pointed matched circle $c$, and so we may identify the parametrized surface $S$ with $S_c$.  

\subsection{Mapping cylinders of parametrized surfaces}

\begin{figure}
\centering
\includegraphics[width=110mm]{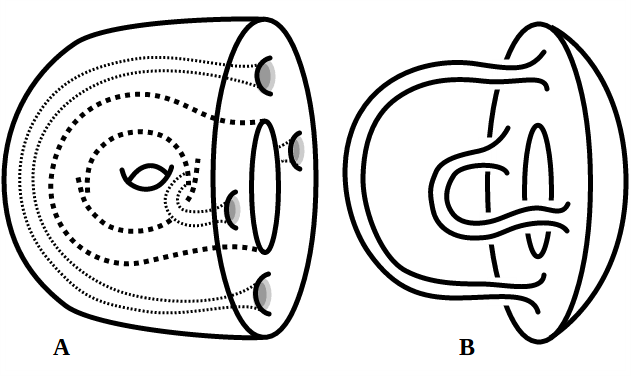}
\caption{A Heegaard decomposition of a mapping cylinder $T^2 \times I$ into two handlebodies, $A$ and $B$.}  
\label{Mapping_Cylinder}
\end{figure}

\begin{figure}
\centering
\includegraphics[width=80mm]{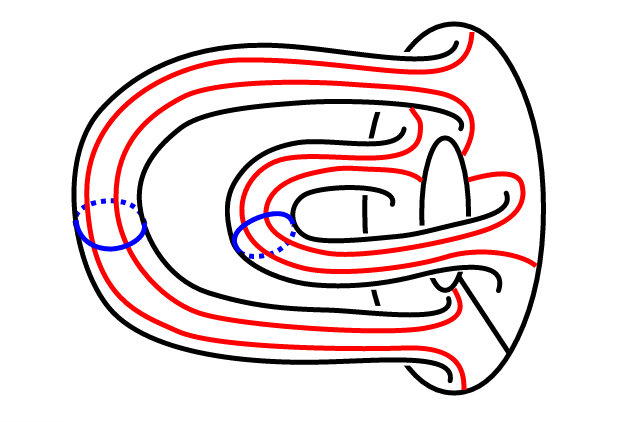}
\caption{Bordered Heegaard diagram for the mapping cylinder associated to the identity function on a parametrized torus.}  
\label{HD_Mapping_Cylinder}
\end{figure}

Let $c$ and $c'$ be pointed matched circles of genus $g$, let $S_c$ and $S_{c'}$ be their associated parametrized surfaces, and let $z$ and $z'$ be the marked points on $c \subset S_c$ and $c' \subset S_{c'}$, respectively.  A mapping cylinder from $c$ to $c'$ is an orientation-preserving diffeomorphism from $S_c$ to $S_{c'}$ which preserves the marked disk and point, where two such diffeomorphisms are considered equivalent if there is an isotopy between them which also preserves the marked disk and point.  \\

Equivalently, a mapping cylinder between $c$ and $c'$ is a class of bordered 3-manifolds diffeomorphic to a Riemann surface cross an interval, with:

\begin{enumerate}

\item
The two boundary components marked ``left" and ``right",

\item
The left boundary component parametrized by $-S_c$ and the right by $S_{c'}$, and

\item
A marked section over the interval, with framing, which includes the two marked points on the boundary, and which extends the framings of $T_z(S_c)$ and $T_{z'}(S_{c'})$ arising from the oriented curves $c$ and $c'$.  

\end{enumerate}

Two such manifolds are equivalent when there is a diffeomorphism between them, taking the left (right) boundary component of one to the left (right) boundary component of the other, which preserves the parametrizations of both boundary components as well as the framed section.  \\

For any genus $g$, we may construct a category in which the objects are the pointed matched circles of genus $g$, and the morphisms are the mapping cylinders.  This category is the mapping class groupoid in genus $g$ \cite{LOT2}.  Note that, given a pointed matched circle $c$, the group of morphisms from $c$ to itself is the mapping class group for the parametrized surface $S_c$.  \\

As we will see, bordered Heegaard Floer homology constructs a functor from the genus $g$ mapping class groupoid to the category of differential graded algebras, with morphisms given by type $DA$ bimodules.  To describe this functor, we will begin by considering bordered Heegaard diagrams associated to mapping cylinders.  \\

To construct a bordered Heegaard diagram for a mapping cylinder, first we choose a parametrization of an interior fiber which is compatible with the marked section.  For a mapping cylinder described as a class of diffeomorphisms $f: S_c \rightarrow S_{c'}$, this means choosing a factorization of $f$ into mapping cylinders $f_l: S_c \rightarrow S, f_r: S \rightarrow S_{c'}$ with $f = f_r \circ f_l$, where $S$ is some parametrized genus $g$ surface.  \\

This allows us to decompose the mapping cylinder into handlebodies $A$ and $B$ as follows.  First, extend the parametrization of the interior fiber to a parametrization of the mapping cylinder by $S \times [0,1]$.  Let $D \subset S$ be the marked disk, and let $T \subset S$ be the genus $g$ subsurface with two boundary components obtained by thickening the marked elements of $S$ and removing the disk.  Then, for some $\epsilon \in (0, 0.5)$, define the handlebodies $A$ and $B$ by:  

\begin{equation}
B = (D \times [0,1]) \cup (T \times [0.5 - \epsilon, 0.5 + \epsilon]), \nonumber
\end{equation}

and

\begin{equation}
A = S \times [0,1] - B. \nonumber
\end{equation}

\begin{figure}
\centering
\includegraphics[width=50mm]{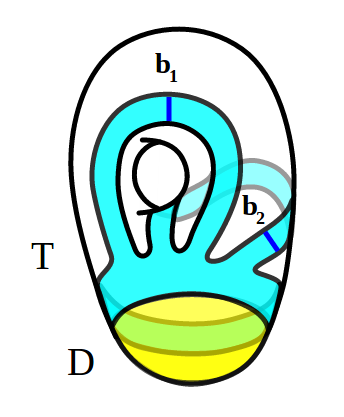}
\caption{A (genus one) interior fiber $S$, with the subsurfaces $T$ and $D$ shaded.  The arcs $b_1$ and $b_2$ are marked, as well.}  
\label{Interior_Fiber}
\end{figure}

Take $\{b_1, \ldots b_g\}$ to be properly embedded arcs in $T$ which separate the thickened arcs in the parametrization of $S$ (see Figure 4).  We may then define the $\beta$ disks by $B_i = b_i \times [0.5 - \epsilon, 0.5 + \epsilon]$.  To construct the $\alpha$ disks, let $a_1, \ldots a_g$ be the arcs in the parametrization of $S_c$, included in $S$ by its identification with $S \times \{0\}$.  We may deform the $a_i$ so that they lie on $T \times \{0\}$, since $T$ is a genus $g$ subsurface of $S$.  The $\alpha$ disks which intersect $S \times \{0\}$ are then given by $A_i = a_i \times [0, 0.5 - \epsilon]$.  Similarly, the $\alpha$ disks which intersect $S \times \{1\}$ are given by $A'_i = a'_i \times [0.5 + \epsilon, 1]$, where the $a'_i$ are the arcs in the parametrization of $S_{c'}$.  \\

After smoothing the corners, we may identify the left half of the Heegaard surface $A \cap B \cap (S \times [0, 0.5])$ with $T$, with arcs $\alpha_1, \ldots \alpha_g$ given by $\alpha_i = f_l(a_i)$.  We can identify the right half of the Heegaard surface with $-T$, with arcs $\alpha'_i = f_r^{-1}(a'_i)$ for $i \in \{1, \ldots g\}$.  The curves $\beta_1, \ldots \beta_g$ are given by $\beta_i = b_i \cup -b_i$.  \\

This construction allows us to build Heegaard diagrams to emphasize any preferred factorization of a mapping cylinder.  In particular, we may take $S = S_c$ and $f_l = I$, in which case the left half of the diagram is standard while the arcs on the right side have been altered by $f^{-1}$, or we may take $S = S_{c'}$ and $f_r = I$ to produce a diagram with a standard right half.  \\

\begin{figure}
\centering
\includegraphics[width=80mm]{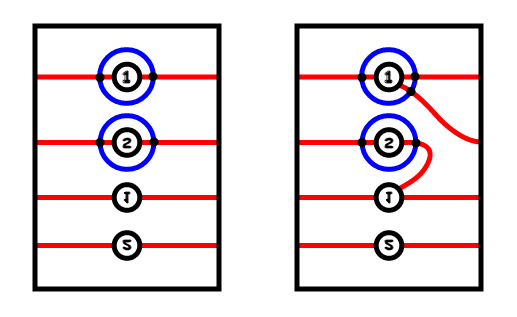}
\caption{Left:  An alternate depiction of the bordered Heegaard diagram from Figure 3.  Right:  A bordered Heegaard diagram representing a nontrivial mapping cylinder of genus one.}  
\label{TwoMappingCylinders}
\end{figure}

Given two Heegaard diagrams for the same mapping cylinder which are constructed from different middle parametrizations, we know we can get from one to the other by a sequence of isotopies, handleslides, stabilizations and destabilizations.  It's useful to look at one method for accomplishing this.  \\

Let $f$ be a mapping cylinder with $f_r \circ f_l = g_r \circ g_l = f$ two factorizations, and let $H$ and $H'$ be the associated Heegaard diagrams, respectively.  To take the $\alpha$ arcs of $H$ to those of $H'$ we apply the diffeomorphism $g_l \circ f_l^{-1}|_T$ to the left half of the diagram, and $g_r^{-1} \circ f_r|_{-T}$ to the right half.  Note that:  

\begin{align*}
g_r^{-1} \circ f_r|_{-T} &=  g_l \circ f^{-1} \circ f_r|_{-T}\\
&= g_l \circ f_l^{-1}|_{-T},
\end{align*}

so we are applying the same diffeomorphism to both halves.  \\

Now consider the following handleslide.  Begin with two arcs in the parametrization of $S$ with a pair of adjacent end points, and let $b_i, b_j$ be the associated arcs in $T$.  The adjacency gives us a curve in the boundary of $T$, running from one end of $b_i$ to one end of $b_j$, which does not intersect any other such end points.  In $H$ this becomes an arc from $\beta_i$ to $\beta_j$, and we may slide $\beta_i$ over $\beta_j$ along this arc.  \\

\begin{figure}
\centering
\includegraphics[width=120mm]{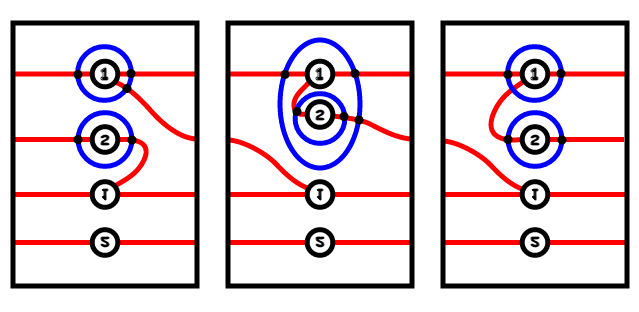}
\caption{Left:  A Heegaard diagram for a mapping cylinder.  Middle:  The same diagram after altering the Heegaard surface by a diffeomorphism.  Right:  A new Heegaard diagram for the same mapping cylinder, obtained by a handleslide.}  
\label{Arcslide}
\end{figure}

This results in a new Heegaard diagram of the form we are using, and it corresponds to altering the parametrization of the interior fiber by an arc slide.  We may do this for any arc slide, and arc slides generate the mapping class groupoid, so we can realize any diffeomorphism in this way.  This allows us to modify the $\alpha$ arcs as desired while keeping the $\beta$ curves in the same form.  

\subsection{The bordered invariants for surfaces and mapping cylinders}

To a pointed matched circle $c$, bordered Heegaard Floer theory associates a differential graded algebra $A_c$ \cite{LOT1}.  If $c$ is a separating curve on a Heegaard diagram $H$ for a manifold $Y$ as in section 3.1, then the invariant $A_c$ contains information about the behavior of $\widehat{HF}(Y)$ near $c$.  Namely, if we have a holomorphic disk in $H \times [0,1] \times \mathbb{R}$, then the restriction of this disk to $c \times [0,1] \times \mathbb{R}$ is a collection of arcs, which we may represent by a strand diagram.  We put additional markings on this diagram to record the behavior of sheets of this disk which do not intersect $c \times [0,1] \times \mathbb{R}$, and the strand diagrams of this form are the generators of $A_c$ over $\mathbb{Z} / 2 \mathbb{Z}$.  \\

In most cases, the product of two strand diagrams is defined as their concatenation if it exists, and $0$ otherwise.  The exception to this is that strand diagrams with double crossings are not permitted, and so if two diagrams have crossings which ``undo" each other, then their product is also defined to be $0$.  The differential of a strand diagram is the sum of all diagrams obtained from resolving one of its crossings, also with the exception that resolutions which undo a second crossing are excluded.  \\

As the strand diagrams represent the behavior of holomorphic disks on the curve $c$, the algebra operations represent the behavior of ends of one-dimensional families of holomorphic disks near this curve.  The proofs that the differential squares to zero and that the operations satisfy the Leibnitz rule arise from counts of the ends of these moduli spaces.  \\

Given a mapping cylinder $f: S_c \rightarrow S_{c'}$, its Heegaard Floer invariant is a type $DA$ bimodule $\widehat{CFDA}(f)$ over $A_c$ and $A_{c'}$ [LOT2].  The type $D$ structure on $\widehat{CFDA}(f)$ is an identification $D: \widehat{CFDA}(f) \rightarrow A_c \otimes \chi$.  Here, $\chi$ is the set of $2g$-tuples of intersection points between the $\beta$ curves and $\alpha$ arcs, where each $\beta$ curve includes exactly one intersection point, and each $\alpha$ arc includes at most one.  \\

For an element $(x, a_1, \ldots a_i) \in \chi \otimes (A_{c'})^{\otimes i}$, the product $m_{i+1}(x, a_1, \ldots a_i)$ arises from counting certain rigid holomorphic surfaces in the manifold $H \times [0,1] \times \mathbb{R}$, where $H$ is the Heegaard surface for $f$.  \\

Given composable mapping cylinders $f$ and $g$, \cite{LOT2} have shown that the product $\widehat{CFDA}(f) \boxtimes \widehat{CFDA}(g)$ is quasi-isomorphic to the bimodule $\widehat{CFDA}(g \circ f)$.  Thus the bordered Heegaard Floer invariants for mapping cylinders of genus $g$ comprise a functor from the mapping class groupoid of genus $g$ to the category of DGA's, with morphisms given by type $DA$ bimodules.  

\section{Cornered Lefschetz fibrations}

\subsection{Cornered Lefschetz fibrations}

\begin{definition}

A cornered Lefschetz fibration, or CLF, is a Lefschetz fibration over the rectangle $[0,1] \times [0,1]$, with a marked, framed section, such that:  

\begin{enumerate}

\item
The vanishing cycles are nonseparating,

\item
The ``bottom edge" (the preimage of $[0,1] \times \{0\}$) and the ``top edge" (the preimage of $ [0,1] \times \{1\}$)
are both identified with mapping cylinders, with the ``left corners" (the fibers over $\{(0,0)\}$ and $\{(0,1)\}$)
identified with the left boundary components, and the ``right corners" (the fibers over $\{(1,0)\}$ and $\{(1,1)\}$)
identified with the right boundary components.

\item
The ``right edge" and ``left edge" are each identified with a parametrized Riemann surface cross interval,

\item
The parametrizations induced by these identifications agree on the corners, and

\item
The framed section over $[0,1] \times [0,1]$ agrees with the framed sections on the edges.  

\end{enumerate}

\end{definition}

Given two cornered Lefschetz fibrations, we consider them equivalent when there is a symplectomorphism between them, which restricts to diffeomorphisms between the respective edges and corners, and which preserves the framed section and the parametrizations of all parametrized fibers. \\

\begin{figure}
\centering
\includegraphics[width=80mm]{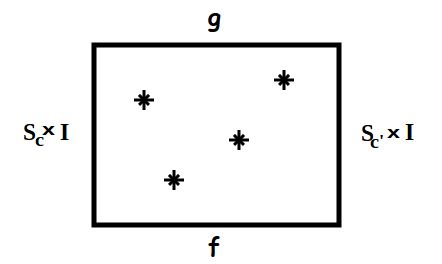}
\caption{A CLF, depicted as its base space $[0,1] \times [0,1]$ with critical values marked.}  
\label{CLF}
\end{figure}

If we restrict our attention to cornered Lefschetz fibrations with a single critical point, we may use an alternate definition.  \\

\begin{definition}

An abstract CLF with one critical point consists of the following data:  

\begin{enumerate}

\item
``Initial" and ``resulting"  abstract mapping cylinders $f, g: S_c \rightarrow S_{c'}$.  

\item
For the initial mapping cylinder, we have a parametrization of an interior fiber given by $f_l$ and $f_r$ with $f_r \circ f_l = f$.

\item
A marked isotopy class of nonseparating simple closed curves $\zeta$ on the parametrized middle fiber.    

\end{enumerate}

This data must satisfy:  

\begin{equation}
g = f_r \circ T_{\zeta} \circ f_l,
\end{equation}

where $T_{\zeta}$ is the negative Dehn twist about $\zeta$, due to our orientation conventions.  

\end{definition}

We consider two such abstract CLF's equivalent if the initial and resulting mapping cylinders are equivalent, and if the identification of the left boundary components of the initial mapping cylinders preserves the preimage of $\zeta$ via $f_l$.  Note that the image of $\zeta$ via $f_r$ is also preserved by the identification of the right boundary components of these mapping cylinders.  

\subsection{Constructing bordered Heegaard triples}

\begin{figure}
\centering
\includegraphics[width=100mm]{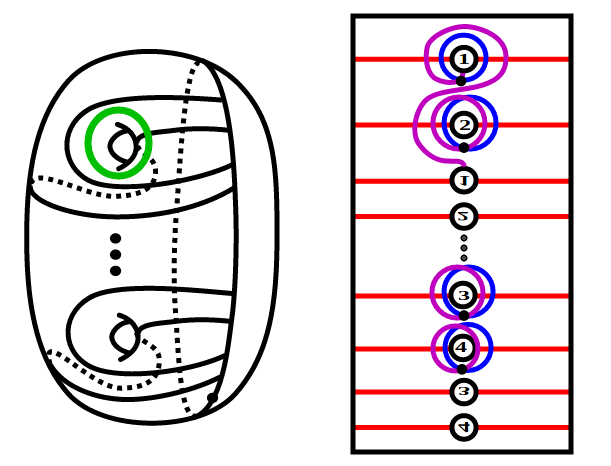}
\caption{Left: The standard parametrized genus g surface $S_g$, with the curve $\eta$ marked in green.  Right:  A bordered Heegaard triple representing the CLF with a single critical point, fiber $S_g$, and vanishing cycle $\eta$.}  
\label{Standard_Curve}
\end{figure}

For a given abstract CLF with one critical point, we may construct a bordered Heegaard triple representing it as follows.  First, choose a factorization so that the curve $\zeta$ is the standard curve $\eta$ on the canonical parametrized surface $S_g$ (see Figure 8).  To see that this is always possible, for a given CLF of this form, $\{f, g, f_l, f_r, \zeta \}$, choose an orientation preserving diffeomorphism $f_l': S_c \rightarrow S_g$, with $f_l^{-1}(\zeta) = f_l'^{-1}(\eta)$.  Now define $f_r': S_g \rightarrow S_{c'}$ by $f_r'= f \circ f_l'^{-1}$.  This new data defines a new abstract CLF by $\{f, g, f_l', f_r', \eta\}$.  Note that $f_r' \circ f_l' = f$, and that:  

\begin{align}
f_r' \circ T_{\eta} \circ f_l' &= f \circ f_l'^{-1} \circ T_{\eta} \circ f_l'  \nonumber \\
&= f \circ T_{f_l'^{-1}(\eta)} \nonumber \\
&= f \circ T_{f_l^{-1}(\zeta)}  \nonumber \\
&= f \circ f_l^{-1} \circ T_{\zeta} \circ f_l \nonumber \\
&= f_r \circ T_{\zeta} \circ f_l  \nonumber \\
&= g. \nonumber
\end{align}

Since $f_l'^{-1}(\eta) = f_l^{-1}(\zeta)$, these CLF's are equivalent. \\

Now we have our CLF expressed as $W = \{f, g, f_l, f_r, \eta\}$.  In order to construct a bordered Heegaard triple representing $W$, start with the diagram for the mapping cylinder $f$, with the middle fiber given by $f_l, f_r$.  By including $\eta$ in this middle fiber, we may interpret it as a knot in the mapping cylinder $f$.  Then $W$ is the cobordism obtained by doing $+1$ surgery on this knot, so we obtain the $\gamma$ curves by altering the $\beta$ curves by a Dehn twist around the projection of $\eta$ to the left half of the Heegaard diagram. \\

\subsection{A morphism associated to abstract CLF's}

Given a bordered Heegaard triple $\{H, \alpha, \beta, \gamma\}$ constructed from the abstract CLF $W = \{f, g, f_l, f_r, \eta\}$, as in the previous section, let $\Theta$ be the tuple of intersection points between the $\beta$ and $\gamma$ curves which generates the highest degree of $\widehat{CFDA}(Y_{\beta, \gamma})$, where $Y_{\beta, \gamma}$ is the 3-manifold with two boundary components obtained from the Heegaard diagram $(H, \beta_i, \gamma_i)$.  Then we have the following definition:  

\begin{definition}

Let $x$ be a generator for $CFDA(Y_{\alpha, \beta})$, and let $y$ be a generator for $CFDA(Y_{\alpha, \gamma})$.  Then a triangle from $x$ to $y$ consists of the following data:  \\

A Riemann surface $S$ with a punctured boundary, along with a proper holomorphic embedding $u: S \rightarrow \bar{H} \times T$.  Here $T$ is a disk with three boundary punctures, with the arcs between the punctures labelled $a, b$ and $c$, and $\bar{H}$ is the completion of the Heegaard surface $H$ obtained by attaching infinite cylindrical ends to the boundary components.  \\

The map $u$ extends continuously to the compactifications of $S$ and $\bar{H} \times T$ obtained by filling the boundary punctures, in a manner which maps the punctures of $S$ to the following points:  

\begin{itemize}

\item
The punctures $x_i \times (a,b)$, where $x_i$ is a point in $x$ and $(a,b)$ is the puncture lying between arcs $a$ and $b$.  

\item
The punctures $y_i \times (a,c)$ and $\Theta_i \times (b,c)$, defined similarly.  

\item
Points of the form $e \times a_t$ or $w \times a_t$, where $a_t$ is some point on $a$, and $e$ and $w$ are the punctures in $\bar{H}$ corresponding to the right and left boundary components of $H$, respectively.  

\end{itemize}

Furthermore, we require that each of the arcs comprising the boundary of $S$ map to a surface of the form $\alpha_i \times a$, $\beta_i \times b$, or $\gamma_i \times c$.  

\end{definition}

With this in mind, we can define a type $DA$ map $F$ associated to our Heegaard triple:  

\begin{definition}

For each generator element $(x, a_1, \ldots a_n) \in \chi \otimes A_{c'}^{\otimes n}$:

$$F_{n+1}(x, a_1, \ldots a_n) = \sum_{y \in \chi', b \in A_c, \Delta \in \Delta_{x, y, b}} b \otimes y.$$

Here $\Delta_{x, y, b}$ is the set of rigid triangles from $x$ to $y$, which approach the Reeb chords $a_1, \ldots a_n$ near $e \times a$ and $b_1 \ldots b_m$ near $w \times a$, such that the product $b_1\ldots b_m = b$.  

\end{definition}

\begin{lemma}
The map $F$ is a morphism of type DA bimodules.  
\end{lemma}

The proof is similar to the proof from \cite{LOT1} that the maps induced by handleslides are chain maps and $A_{\infty}$ maps.  The proof in question involves identifying ends of one-dimensional moduli spaces of triangles, but is complicated by the appearance of triangles with corners at Reeb chords within these ends.  For our purposes this is not an issue, since there are no $\beta$ arcs or $\gamma$ arcs, and so triangles of this type do not exist.  \\

Given two bordered Heegaard triples $H$ and $H'$ for equivalent CLF's, constructed as described above, we may obtain the $\alpha$ arcs of $H'$ from those of $H$ by applying a diffeomorphism to one side of the diagram and its inverse to the other side. To preserve the $\beta$ and $\gamma$ curves as well, we can realize this diffeomorphism by a sequence of handleslides.  Since the diagrams are equivalent the diffeomorphism fixes the projection of $\eta$, and so we may perform these handleslides away from the curves $\beta_1$ and $\gamma_1$.  

\begin{lemma}
Consider a bordered Heegaard triple in which $\beta_1$ and $\gamma_1$ differ by a Dehn twist, and $\beta_i$ and $\gamma_i$ differ by a Hamiltonian isotopy for each $i \neq 1$.  If we perform a sequence of simultaneous handleslides among the $\beta_i$ and $\gamma_i$ for $i \neq 1$, then this will not alter the homotopy class of the induced map.  
\end{lemma}

To prove this, assume the Heegaard triples $H$ and $H'$ differ by a single handleslide.  We must show that the morphisms $F' \circ h$ and $h' \circ F$ are chain homotopic, where $F$ is the triangle map induced by the diagram $H$, $F'$ is the map induced by $H'$, and $h$ and $h'$ are the quasi-isomorphisms induced by the handleslides in question.  \\

The argument is similar to the proof of handleslide invariance for the cobordism map in \cite{OS2}.  First, construct a Heegaard quadruple $\{H'', \alpha, \beta, \gamma, \gamma'\}$ where $\{H'', \alpha, \beta, \gamma\}$ is the triple diagram $H$, and the $\gamma'$ curves are obtained from altering the $\gamma$ curves by the relevant handleslide.  We may compose the triangle maps induced by the diagrams $H$ and $\{H'', \alpha, \gamma, \gamma'\}$.  However, there is an associativity result for such maps, which shows that this is homotopic to the composition of maps induced by the diagrams $\{H'', \alpha, \beta, \gamma'\}$ and $\{H'', \beta, \gamma, \gamma'\}$.  \\

More precisely, we may consider holomorphic curves in $Q \times \bar{H}$, where $Q$ is a disk with four boundary punctures, with the arcs between them labelled $a$, $b$, $c$ and $c'$, and corresponding boundary conditions $\alpha_i \times a$, $\beta_i \times b$, $\gamma \times c$, and $\gamma' \times c'$.  By counting rigid curves of this form, we may define a chain homotopy between the two compositions described above.  The fact that this map is such a chain homotopy arises from counts of the ends of one-dimensional moduli spaces of curves of this type.  Degenerations into two triangles correspond to terms in a composition, and degenerations into quadrilaterals and disks correspond to terms from the map in question followed by or preceded by a differential.  \\

The map induced by the Heegaard triple $\{H'', \beta, \gamma, \gamma'\}$ takes the generators $\Theta_{\beta, \gamma}$ and $\Theta_{\gamma, \gamma'}$ to the generator $\Theta_{\beta, \gamma'}$, and so the composition $h' \circ F$ is homotopic to the map induced by $\{H'', \alpha, \beta, \gamma'\}$.  A similar argument shows that $F' \circ h$ is homotopic to this map as well.  \\

We also have the following result:  

\begin{lemma}
Suppose we have a bordered Heegaard triple as in Lemma 4, and that we slide an $\alpha$ arc or curve over an $\alpha$ curve.  Then this will not change the homotopy class of the induced map.  
\end{lemma}

The argument is similar to the proof of Lemma 4.6, however the associativity result for triangle maps has an additional complication.  This stems from the fact that we are considering a Heegaard quadruple $\alpha, \alpha', \beta, \gamma$ in which the first two sets of curves both interact with the boundary.  As before we define a chain homotopy by counting rigid quadrilaterals with appropriate boundary conditions, and we prove that this map is the desired chain homotopy by counting degenerate quadrilaterals.  However, these degenerate curves may now include punctures which map to points of the form $e \times (a, a')$ or $w \times (a, a')$, where $(a, a')$ is the puncture on the boundary of $Q$ which typically maps to $\Theta_{\alpha, \alpha'} \times (a, a')$.  \cite{LOT1} demonstrated that curves of this type do not contribute to the map, and so the result follows.  \\

\begin{lemma}
Given a bordered Heegaard triple as constructed above, the induced map is independent of the chosen almost-complex structures, and invariant under isotopies of the Heegaard diagram.  
\end{lemma}

In order to prove invariance with respect to the choice of almost-complex structure, we construct a homotopy between the moduli spaces for different almost-complex structures.  This is very similar to Proposition 6.16 of \cite{LOT2} (see also sections 6.4 and 7.4 of \cite{LOT1}).  Given two almost-complex structures $J_0$ and $J_1$ with a one-dimensional family of almost-complex structures $J_t$ between them, there are quasi-isomorphisms $\Phi_{J_0, J_1}$ between the appropriate bimodules.  These maps come from counts of index 0 holomorphic curves in $\bar{H} \times [0,1] \times \mathbb{R}$, in which the almost-complex structure varies with the coordinate in $\mathbb{R}$ and interpolates from $J_0$ to $J_1$.  \\

Denoting by $F_{J_0}$ and $F_{J_1}$ the triangle maps induced by the Heegaard triple for different complex structures, we need to show that $\Phi_{J_0, J_1} \circ F_{J_0}$ is homotopic to $F_{J_1} \circ \Phi_{J_0, J_1}$.  To construct a chain homotopy between these maps, we consider holomorphic maps to $\bar{H} \times T$, where the almost-complex structure depends on the point in $T$, and agrees with $J_0$ near the punctures $(a, b)$ and $(b, c)$ and with $J_1$ near $(a, c)$.  We may then allow this almost-complex structure to vary in a one-parameter family, interpolating between the product complex structure determined by $J_0$ and that determined by $J_1$.  By counting the ends of the resulting parametrized moduli spaces, we can verify that the map in question is the desired chain homotopy.  \\

The argument for invariance with respect to Hamiltonian isotopies is similar.  

\section{A cobordism map and invariance}

\subsection{Horizontal and vertical composition}

\begin{figure}
\centering
\includegraphics[width=60mm]{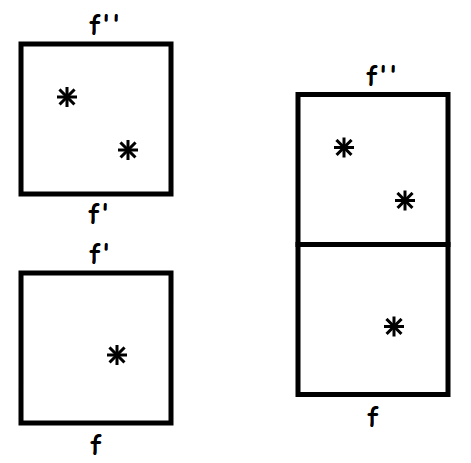}
\caption{Vertical composition of CLF's.}  
\label{Vertical}
\end{figure}

Cornered Lefschetz fibrations may be composed both horizontally and vertically.  Given two CLF's $W$ and $W'$, if the resulting mapping cylinder of $W$ is equivalent to the initial mapping cylinder of $W'$, then there is a unique CLF obtained by identifying $W$ and $W'$ along that mapping cylinder. This is the vertical composition of $W$ and $W'$, written $W \circ_v W'$.  If $Z$ and $Z'$ are CLF's and the fibers in the right edge of $Z$ and the left edge of $Z'$ are parametrized by the same pointed matched circle, then we may identify those edges to define the horizontal composition $Z \circ_h Z'$.  \\

\begin{figure}
\centering
\includegraphics[width=100mm]{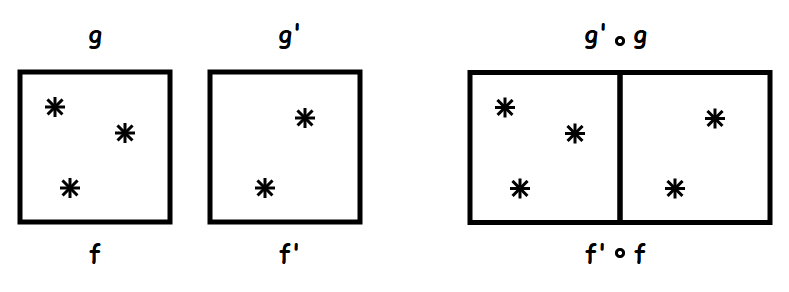}
\caption{Horizontal composition of CLF's.}  
\label{Horizontal}
\end{figure}

In the first case, if we have type $DA$ bimodule maps $F$ and $F'$ associated to $W$ and $W'$ respectively, then we may associate the map $F' \circ F$ to the vertical composition of $W$ and $W'$.  In the case of horizontal composition, suppose $Z$ and $Z'$ have initial mapping cylinders $f$ and $f'$ and resulting mapping cylinders $g$ and $g'$.  If we have type DA maps $G$ and $G'$ associated to $Z$ and $Z'$, then there is an induced map on the tensor product:

\begin{equation}
G \boxtimes G' : CFDA(f) \boxtimes CFDA(f') \rightarrow CFDA (g) \boxtimes CFDA(g').  \nonumber
\end{equation}

Since $CFDA(f) \boxtimes CFDA(f')$ is quasi-isomorphic to $CFDA(f' \circ f)$, and since $CFDA (g) \boxtimes CFDA(g')$ is quasi-isomorphic to $CFDA(g' \circ g)$, we may associate the map $G \boxtimes G'$ to the horizontal composition of $Z$ and $Z'$.  \\

Given a CLF $W$ with initial and resulting mapping cylinders $f$ and $g$, we may express $W$ as a sequence of horizontal and vertical compositions of CLF's, each with at most one critical point.  Such a decomposition of $W$ induces a type DA map $F: CFDA(f) \rightarrow CFDA(g)$.  In the rest of this section we will prove the following result:  

\begin{theorem}
The homotopy class of the map $F$ depends only on the symplectic structure of the CLF $W$.  
\end{theorem}

\subsection{Invariance for CLF's with a single critical point}

First, observe that this result holds for CLF's with no critical points.  This follows from Lemma 1.  \\

Now let $W$ be a CLF with one critical point, expressed as $\{f, g, S, \eta, f_l, f_r\}$, with induced map $F: CFDA(f) \rightarrow CFDA(g)$.  If we express this CLF as a vertical composition then the new induced map will be either $F \circ I_g$ or $I_f \circ F$, both of which are homotopic to $F$, and so we will consider a horizontal decomposition $W = U \circ_h V$.  \\

First, we will assume that $U$ contains a critical point and that $V$ is trivial.  Then these CLF's are of the form $U = \{(f_r^2)^{-1} \circ f, (f_r^2)^{-1} \circ g, S, \eta, f_l, f_r^1\}$ and $V = I(f_r^2)$, for some factorization $f_r = f_r^2 \circ f_r^1$.  Let us further assume that $f^2_r = f_r$, giving us $U = \{f_l, (f_r)^{-1} \circ g, S, \eta, f_l, I \}$ and $V = I(f_r)$.  This induces a type $DA$ map:  

$$F': CFDA(f_l) \boxtimes CFDA(f_r) \rightarrow CFDA((f_r)^{-1} \circ g) \boxtimes CFDA(f_r).$$

These bimodules are quasi-isomorphic to $CFDA(f)$ and $CFDA(g)$, respectively, and we would like to show that the maps $F$ and $F'$ are homotopic.  \\

\begin{figure}
\centering
\includegraphics[width=80mm]{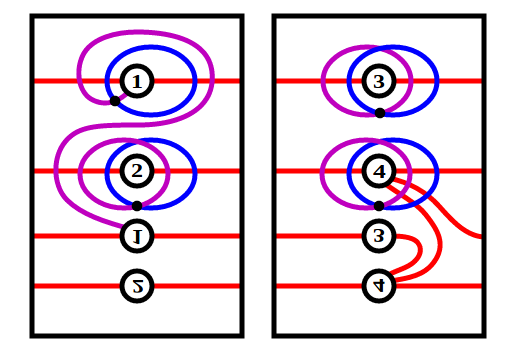}
\caption{The Heegaard triples $H_1$ (left) and $H_2$ (right).}  
\label{Pairing}
\end{figure}

Let $H$ be the Heegaard triple for $W$ arising from its description.  Let $H_1$ be the Heegaard triple for $U$ obtained from the description $U = \{f_l, (f_r)^{-1} \circ g, S, \eta, f_l, I \}$, and let $H_2$ be the Heegaard triple for $V$ defined by the factorization $f_r = f_r \circ I$.  Construct a new Heegaard triple $H'$ by identifying the right boundary component of $H_1$ with the left boundary component of $H_2$.  This is a Heegaard triple which represents $W$, although its genus is higher than that of $H$.  \\

Let $F''$ be the type $DA$ bimodule map induced by the Heegaard triple $H'$.  Then we have the following lemma:  

\begin{lemma}
(Stabilization)  The maps $F$ and $F''$ are chain homotopic.  
\end{lemma}

\begin{figure}
\centering
\includegraphics[width=120mm]{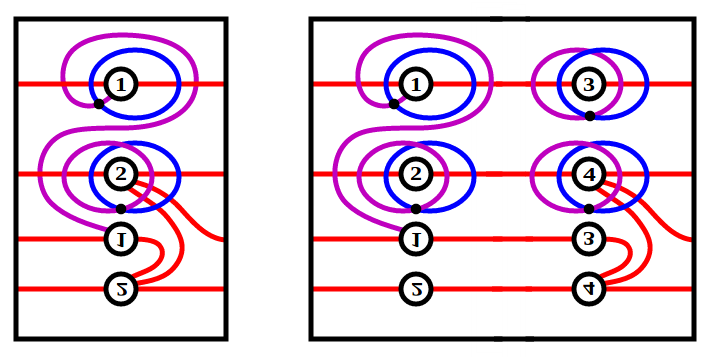}
\caption{The Heegaard triples $H$ (left) and $H'$ (right).}  
\label{Stabilization}
\end{figure}

To show this we will obtain the Heegaard triple $H$ from $H'$ by a certain sequence of handleslides and destabilizations, and show that these moves do not change the homotopy class of the induced map.  The diagram $H'$ has $2g$ $\alpha$ curves, along with $4g$ $\alpha$ arcs, and each $\alpha$ curve intersects two $\beta$ curves and two $\gamma$ curves once.  Call these curves $\alpha_i, \beta_i, \beta_i', \gamma_i, \gamma_i'$ for each $i \leq 2g$.  \\

For each $i$, let $\beta_i'$ be the curve which intersects $\alpha$ arcs with end points on the right boundary component of $H'$, and let $\gamma_i'$ be the analogous $\gamma$ curve.  We may remove these intersections by sliding the arcs over $\alpha_i$, along a segment of $\beta_i'$.  Next, we slide $\beta_i$ over $\beta_i'$ along a segment of $\alpha_i$, while simultaneously sliding $\gamma_i$ over $\gamma_i'$ along the analogous arc.  The proof that this move does not change the homotopy class of the map is similar to the proof of Lemma 4.  \\

Following these handleslides, for each $i$ the curve $\alpha_i$ intersects $\beta_i'$ and $\gamma_i'$ once, and $\beta_i'$ and $\gamma_i'$ differ by a Hamiltonian isotopy, but this triple is disjoint from all other curves.  We wish to destabilize the diagram by removing each such triple.  We may do this if $\alpha_i, \beta_i'$ and $\gamma_i'$ lie in the region of the diagram containing the marked arc, since there is a one-to-one correspondence between generators, rigid disks, and rigid triangles before and after such a destabilization, in this case.  \\

For a triple $\alpha_i, \beta_i', \gamma_i'$, there is a path from an intersection of $\beta_i'$ and $\gamma_i'$ to the marked arc, which does not intersect $\beta_1$ or $\gamma_1$.  This path may cross $\alpha$ curves or arcs, or other $\beta$ curves and their analogous $\gamma$ curves.  If the first crossing is with a $\beta$ and $\gamma$ curve, we may remove it by sliding these curves over $\beta_i'$ and $\gamma_i'$ along the path, and then sliding them again over $\beta_i'$  and $\gamma_i'$ along $\alpha_i$ to remove the intersection created by the previous slide.  If the first crossing is with a curve or arc $\alpha_j$, we may deform $\beta_i'$ and $\gamma_i'$ by a finger move along the path so that they each intersect $\alpha_j$ twice, and then remove these intersections by sliding $\alpha_j$ over $\alpha_i$ twice, along the two segments of $\beta_i'$ which join them.  \\

A sequence of moves of this type will bring the triple $\alpha_i, \beta_i', \gamma_i'$ to the region adjacent to the marked arc, while leaving them disjoint from all other curves, and so we may then destabilize the diagram without changing the homotopy class of the induced map.  Since the previous moves were all handleslides over $\alpha_i, \beta_i'$ or $\gamma_i'$, the resulting destabilized diagram is isotopic to the diagram obtained by removing the triple without performing these handleslides.  Thus we may perform this destabilization for each $i \leq 2g$, obtaining the diagram $H$ without changing the homotopy class of the resulting map.  \\

Now we need the following result:  

\begin{lemma}

(Pairing)  The maps $F''$ and $F'$ are chain homotopic.  

\end{lemma}

First note that $F'$ is the tensor product $G \boxtimes I(f_r)$, where $G$ is obtained by counting triangles on the diagram $H_1$.  The identity map I has no higher maps, and so the higher maps of $G \boxtimes I$ are of the form:  

$$(G \boxtimes I)_{i+1}(x \otimes y, a_1, \ldots a_i) = \sum  F_{j+1}(x, b_1, \ldots b_j) \otimes z,$$

where $\{ b_1, \ldots b_j, z \}$ is a term arising in the type $D$ product of $y$ with $\{ a_1, \ldots a_i \}$.  \\

These terms correspond to counts of rigid triangles in the Heegaard triple $H_1$, and rigid disks in the Heegaard diagram $\{H_2, \alpha, \beta \}$.  Specifically, suppose the expression $\sum  F_{j+1}(x, b_1, \ldots b_j) \otimes z$ includes the term $\{ c_1, \ldots c_k, u \} \otimes v$.  Then there are an odd number of collections of rigid triangles in $H_1$ and rigid disks in $\{H_2, \alpha, \beta \}$ which represent this term and are compatible.  \\

For each rigid disk in $\{H_2, \alpha, \beta \}$ there is a family of triangles in $H_2$, obtained by replacing each $\beta$ edge with the analogous concave corner between $\beta$ and $\gamma$.  In the degenerate limit where the $\beta$ and $\gamma$ curves of $H_2$ strictly coincide, this family of triangles would be obtained by switching from the $\beta$ curve to the corresponding $\gamma$ curve at any time $t \in \mathbb{R}$ along the
$\beta$-edge. The actual family of triangles we consider is obtained by deforming these via a Hamiltonian isotopy of the $\gamma$ curves. On the given Heegaard triple, this means that at the chosen point along the $\beta$-edge we jump from the $\beta$-curve to the $\gamma$-curve, by attaching a thin triangle ending at the intersection point
$\Theta_{\beta,\gamma}$. The resulting degrees of freedom yield an odd number of rigid triangles whose west degenerations occur at the appropriate time.  We may glue these triangles to the triangles in $H_1$, thus obtaining an odd number of rigid triangles in the destabilized diagram.  \\

Conversely, suppose we have such a rigid triangle.  Its domain is a union of triangles in the diagrams for $H_1$ and for $H_2$.  Since each nontrivial triangle in $H_2$ corresponds to a disk in $\{H_2, \alpha, \beta \}$, these triangles all represent families of dimension greater than or equal to one.  Therefore the corresponding triangles in $H_1$ must be rigid.  A count of the dimensions of the triangles in $H_2$ shows that the analogous disks must be rigid as well.  \\

Now we may prove the following:  

\begin{lemma}

Given a decomposition of a CLF with a single critical point, the homotopy class of the induced type DA map does not depend on the decomposition.  

\end{lemma}

We may relax our initial assumptions, and allow the mapping cylinders $f_r$ and $f_r^2$ to differ.  We may decompose $U$ as $U = U_1 \circ_h U_2$, where $U_1 = \{f_l, (f_r)^{-1} \circ g, S, \eta, f_l, I\}$ and $U_2 = I(f_r^1)$, and then express $W$ as $U_1 \circ_h U_2 \circ_h V$, where $V=I(f_r^2)$.  This decomposition satisfies our previous assumptions, and so the map induced by the decomposition is homotopic to $F'$.  The invariance result for CLF's with no critical points, along with Lemma 1, show that the map induced by $U_1 \circ_h I(f_r)$ is homotopic to $F'$ as well.  This decomposition also satisfies our initial assumptions, and so $F$ and $F'$ are homotopic.  \\

The case of a horizontal decomposition $W = U \circ_h V$ where $U$ is trivial and $V$ has a single critical point is similar; while the formula for $I \boxtimes G$ is different, the underlying geometric argument is essentially the same.  \\

\subsection{Interchanging horizontal and vertical compositions}

\begin{lemma}

(Horizontal versus vertical)  Let $W$ be a CLF with at least two critical points, expressed as a composition of CLF's each with a single critical point, and let $F$ be the type $DA$ map induced by this decomposition.  Then there is a purely horizontal decomposition of $W$ which induces the same map up to homotopy.  

\end{lemma}

Proof:  It suffices to show that any individual vertical composition may be removed or replaced with a horizontal composition, without altering the homotopy class of the resulting type $DA$ map.  With that in mind, assume that $W$ is expressed as a vertical composition $W_1 \circ_v W_2$, where the cLf $W_1$ has initial mapping cylinder $f$ and resulting mapping cylinder $f'$, and $W_2$ has initial mapping cylinder $f'$ and resulting mapping cylinder $f''$.  \\

\begin{figure}
\centering
\includegraphics[width=120mm]{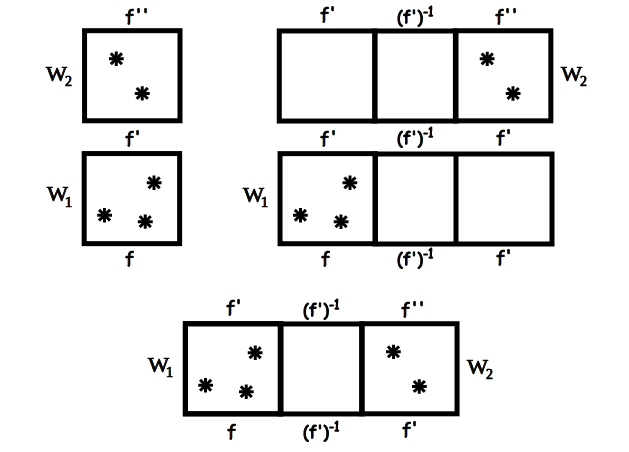}
\caption{Illustration of the steps involved in proving Lemma 5.5.}  
\label{HVV}
\end{figure}

We will argue that the composition $W_1 \circ_h I_{(f')^{-1}} \circ_h W_2$, which induces a type $DA$ map from $f' \circ (f')^{-1} \circ f = f$ to $f'' \circ (f')^{-1} \circ f' = f''$, yields the same map as $W_1 \circ_v W_2$ up to homotopy (See Figure 13).  First, note that $W_1$ and $W_1 \circ_h I_{Id}$ induce homotopic maps, as do $W_2$ and $I_{Id} \circ_h W_2$.  The latter may be decomposed as $I_{f'} \circ I_{(f')^{-1}} \circ W_2$, and by the invariance result for cLf's with no critical points, this change does not alter the homotopy class of the induced map.  \\

Next, we may apply Lemma 2.2 to show that the map induced by $$(W_1 \circ_h I_{Id}) \circ_v (I_{f'} \circ I_{(f')^{-1}} \circ W_2)$$ is homotopic to the map induced by $$(W_1 \circ_v I_{f'}) \circ_h (I_{Id} \circ_v (I_{(f')^{-1}} \circ W_2).$$

By Lemma 2.1, this map is homotopic to that induced by $W_1 \circ_h I_{(f')^{-1}} \circ_h W_2$, as desired.  

\subsection{Invariance under Hurwitz moves}

We have now demonstrated that any decomposition of a CLF may be replaced with a purely horizontal decomposition, without altering the homotopy class of the induced map.  It remains to show that any two horizontal decompositions of the same CLF induce homotopic maps. \\

Given such a horizontal decomposition, there is an ordering of the critical points from ``left" to ``right", according to where they occur in the decomposition.  If two horizontal decompositions of the same CLF result in the same ordering of critical points, then we may construct a common refinement of these decompositions.  Since horizontal compositions of type $DA$ maps are associative up to homotopy, we may use the invariance result for CLF's with a single critical point to show that these two compositions induce homotopic maps.  \\

Now we will show that two horizontal decompositions of the same CLF induce the same map up to homotopy, even if they order the critical points differently.  It is sufficient to treat the case in which these orderings differ by a transposition.  Let $W$ and $W'$ be two CLF's each with a single critical point, which may be composed horizontally as $W \circ_h W'$.  We may express $W$ as an abstract CLF with $f_r = Id$, and $W'$ as an abstract CLF with $f_l = Id$.  \\

\begin{figure}
\centering
\includegraphics[width=120mm]{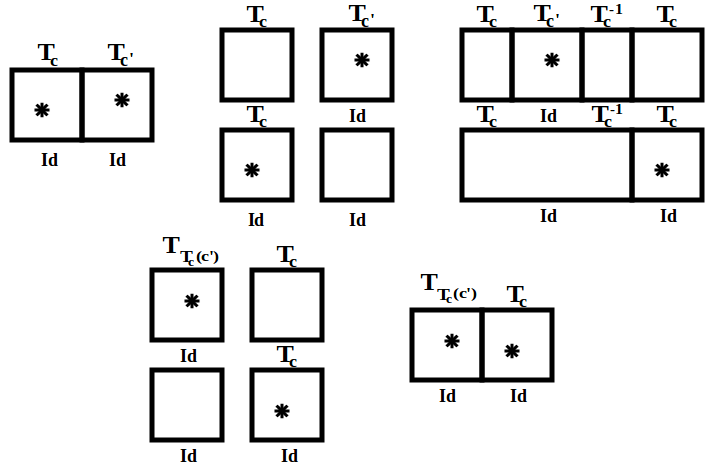}
\caption{Illustration of the steps involved in proving Lemma 5.6.}  
\label{Hurwitz_move}
\end{figure}

If $f$ and $f'$ are the initial mapping cylinders of $W$ and $W'$, respectively, then we may decompose $W$ as $I_{f} \circ_h V$ and $W'$ as $V' \circ_h I_{f'}$, where $V$ and $V'$ are CLF's each with a single critical point, both from the identity to a Dehn twist.  It then suffices to show the following:   

\begin{lemma}

(Hurwitz move)  There is an alternate horizontal decomposition of $V \circ_h V'$, which induces the same map up to homotopy, and which reverses the ordering of the two critical points.  

\end{lemma}

Let $T_{c}$ and $T_{c'}$ be the resulting mapping cylinders of $V$ and $V'$, respectively.  Then we may decompose $V$ as $V \circ_v I_{T_{c}}$, and $V'$ as $I_{Id} \circ_v V'$.  By applying Lemma 2.2, we can then show that $V \circ_h V'$ induces the same map, up to homotopy, as $V \circ_v (I_{T_{c}} \circ_h V')$.  \\

By Lemma 2.1, this map is homotopic to the map induced by:  

$$(I_{Id} \circ_h V) \circ_v (I_{T_{c}} \circ_h V' \circ_h I_{(T_{c})^{-1}} \circ_h I_{T_{c}}).$$

However, the CLF $I_{T_{c}} \circ_h V' \circ_h I_{(T_{c})^{-1}}$ is equivalent to a CLF $V''$, with one critical point, from the identity function to the Dehn twist $T_{T_{c}(c')}$.  By another application of Lemma 2.2, the map induced by $V \circ_h V'$ is thus homotopy equivalent to the map induced by $V'' \circ_h V$.  These two CLF's differ by a Hurwitz move, which preserves the symplectic structure but reverses the order of the two critical points.  This completes the proof of Theorem 5.1.  

\section{Applications}

\subsection{The invariant as a 2-functor}

For each genus $g$ the mapping class groupoid of genus $g$ may be extended to a 2-category, by taking cornered Lefschetz fibrations to be the 2-morphisms.  We may also consider the 2-category whose objects are differential graded algebras, with 1-morphisms given by quasi-isomorphism classes of type $DA$ bimodules, and 2-morphisms given by chain homotopy classes of type $DA$ morphisms.  With this in mind, we have the following theorem:  

\begin{theorem}

The bordered invariants for surfaces and mapping cylinders, along with the maps induced by CLF's, comprise a 2-functor.  

\end{theorem}

This is almost directly a consequence of the invariance result, as we will see.  \\

Recall that, given 2-categories $C$ and $C'$, a 2-functor $F: C \rightarrow C'$ consists of the following data:  

\begin{enumerate}

\item
For each object $x$ in $C$, an object $F(x)$ in $C'$.  

\item
For each morphism $f: x \rightarrow y$ in $C$, a morphism $F(f): F(x) \rightarrow F(y)$ in $C'$.  

\item
For each 2-morphism $\phi: f \rightarrow g$ in $C$, a 2-morphism $F(\phi): F(f) \rightarrow F(g)$ in $C'$.  

\end{enumerate}

This data must satisfy:  

\begin{enumerate}

\item
$F$ preserves identity morphisms and 2-morphisms.  This means that for every object $x$ in $C$ we have $F(I_x) = I_{F(x)}$, and for every morphism $f$ in $C$ we have $F(I_f) = I_{F(f)}$.  

\item
$F$ preserves composition of morphisms, so $F(f_1 \circ f_2) = F(f_1) \circ F(f_2)$, for any composable morphisms $f_1$ and $f_2$ in $C$.  

\item
$F$ preserves both horizontal and vertical composition of 2-morphisms.  This means that, given morphisms $f_1, g_1 : x \rightarrow y$ and $f_2, g_2 : y \rightarrow z$ in $C$, and 2-morphisms $\phi_i : f_i \rightarrow g_i$, we have that $F(\phi_2 \circ_h \phi_1) = F(\phi_2) \circ_h F(\phi_1)$.  Furthermore, given a morphism $h_1: x \rightarrow y$ and a 2-morphism $\psi_1: g_1 \rightarrow h_1$, we also have that $F(\psi_1 \circ_v \phi_1) = F(\psi_1) \circ_v F(\phi_1).$

\end{enumerate}

In our case, the 2-functor $F$ takes a parametrized surface $\Sigma$ to the DGA $A_{\Sigma}$, a mapping cylinder $f: \Sigma_1 \rightarrow \Sigma_2$ to the type $DA$ bimodule $CFDA(f)$ over $A_{\Sigma_1}$ and $A_{\Sigma_2}$, and a CLF $W$ between mapping cylinders $f$ and $g$ to the induced map $F(W): CFDA(f) \rightarrow CFDA(g)$.  All of these associations are up to quasi-isomorphism and chain homotopy, and so $F$ is well-defined.  \\

\cite{LOT2} demonstrated that the bimodule $CFDA(I_{\Sigma})$ over two copies of $A_{\Sigma}$ is quasi-isomorphic to $A_{\Sigma}$ as a type $DA$ bimodule over itself.  We have seen that the map induced by a trivial Heegaard triple for a CLF with no critical points is equal to the identity map on the appropriate bimodule, and so we can see that $F$ satisfies the first criterion.  \cite{LOT2} have also shown that the bimodules $CFDA(f) \boxtimes CFDA(g)$ and $CFDA(g \circ f)$ are quasi-isomorphic, and so $F$ meets the second criterion as well.  \\

To see that $F$ preserves both types of composition of 2-morphisms, note that we defined $F(W)$ to be the map induced by any horizontal or vertical decomposition of $W$, and then showed that the choice of decomposition doesn't matter.  This demonstrates that $F$ is a 2-functor, proving Theorem 6.1.  

\subsection{Calculating the invariant}

Given a CLF $W$ with $n$ critical points and fibers of genus $g$, we may express $W$ as a horizontal composition of the following form:  

$$I(f_1) \circ_h W_g \circ_h I(f_2) \circ_h \ldots \circ_h W_g \circ_h I(f_{n+1}),$$

where $W_g$ is any given CLF with a single critical point and genus $g$ fibers.  \\

This shows that, in order to calculate the map associated to any CLF with fibers of genus $g$, it suffices to know the bimodules associated to mapping cylinders of that genus, and the map associated to a single CLF $W_g$.  \cite{LOT3} have shown that we may calculate $CFDA(f)$ for any mapping cylinder $f$ provided that we have a decomposition of $f$ into arc slides.  Thus the calculation of the map associated to a single CLF with one critical point in each genus would provide the remaining necessary piece.  

\section{Further remarks}

Broken fibrations are a natural generalization of Lefschetz fibrations, in which we allow for smooth one-dimensional families of singular fibers, as well as the usual isolated singular fibers, and in which the genus of the fibers difers by one on either side of such a family.  While Lefschetz fibrations are necessarily symplectic, any smooth 4-manifold may be represented by a broken fibration \cite{AK, Baykur, GK, Lekili}.  By defining cobordism maps associated to broken fibrations, it should be possible to generalize the results of this thesis to obtain a full 2+1+1 TQFT.  \\

This problem is tractable because broken fibrations, like Lefschetz fibrations, may be decomposed into elementary pieces.  One of these pieces is a trivial cobordism between a certain 3-manifold $Y$ and itself.  Here $Y$ is any cobordism between a parametrized genus $g$ surface and a parametrized genus $g+1$ or $g-1$ surface, provided that $Y$ arises from adding a one-handle or two-handle, respectively.  \\

The other new elementary pieces are 4-manifolds with corners that come from adding one-handles and three-handles.  The appropriate cobordism maps for such pieces are analogous to the maps \cite{OS2} developed for one-handle and three-handle additions between closed 3-manifolds.  \\

Once these components are in place, one can attempt to prove that the resulting maps associated to general cobordisms with corners do not depend on the choice of decomposition.  Lekili \cite{Lekili} developed a collection of moves for modifying broken fibrations without altering their smooth structures, and Williams \cite{Williams} proved that these moves are sufficient to relate any two mutually homotopic broken fibrations which represent the same 4-manifold.  It would be desirable to study the behavior of the cobordism maps as we apply these moves, with the hope that the resulting maps will be homotopic.

\end{document}